\documentclass[]{article}

\usepackage{amsthm} % For the `proof' environment

\usepackage{underoverlap} % For contractions / expansions

\usepackage{latexsym}
\usepackage{xypic}
\usepackage{amsmath}
\usepackage{amssymb}
\xyoption{all} % for odd diagrams ...

\newcounter{thm}

\newtheorem{theorem}[thm]{Theorem}

\newtheorem{lemma}[thm]{Lemma}

\newtheorem{definition}[thm]{Definition}

\newtheorem{proposition}[thm]{Proposition}

\newtheorem{corollary}[thm]{Corollary}

\newtheorem{remark}[thm]{Remark}

\newtheorem{notation}[thm]{Notation}

\newcommand{\IN}{{\mathcal I (\mathbb N)}}
\newcommand{\Fk}{{\mathcal F _k}}
\newcommand{ \monus }{\stackrel{\bullet}{-}}
\newcommand{ \KBN }{{ \mathcal N_{k^{(\_)}}}}

\newcommand{\lcm}{{\text{\upshape lcm}}}
\newcommand{\md}{{\text{\upshape mod}}}

\newcommand{\pr}{{\mbox{\tt primes}}}

\begin{document}

\author{{Peter M. Hines}}% \\ \ \\ Peter  Hines}
% \title{``$p$'' is for ``{\em polycyclic}''  --  \\ inverse monoids in number theory}
 \title{The Inverse Semigroup Theory \\ of Elementary Arithmetic}%Une g\'{e}n\'{e}ralisation des mono\"{i}des polycycliques}

\maketitle

\abstract{We curry the elementary arithmetic operations of addition and multiplication to give monotone injections on $\mathbb N$, and describe \& study the inverse monoids that arise from also considering their generalised inverses. This leads to well-known classic inverse monoids, as well as a novel inverse monoid (the `arithmetic inverse monoid' $\mathcal A$) that generalises these in a natural number-theoretic manner.
	
Based on this, we interpret classic inverse semigroup theoretic concepts arithmetically, and vice versa. Composition and normal forms within $\mathcal A$ are based on the Chinese remainder theorem, and a minimal generating set corresponds to all prime-order polycyclic monoids.  This then gives a close connection between Nivat \& Perot's normal forms for polycyclic monoids, mixed-radix counting systems, and $p$-adic norms \& distances.
}
\section{Introduction}
A simple but remarkably powerful notion in logic and theoretical computer science is that of {\em currying} : replacing a multi-argument function by an indexed family of single-argument functions. Given some data-type $D$, and a function of type $f:D\times D\rightarrow D$, we replace the single multi-argument function by a $D$-indexed family of single-argument functions $\{ f_a:D\rightarrow D\}$, defined by $f_a(\_ )=f(a,\_)$, for all $a$ of type $D$.

The name `currying' comes from Haskell Curry, who used this simple technique to great effect in his work on lambda calculus.  As well as its use in computer science, it forms a core part the theory of Cartesian closed categories \cite{LS}, and under the Curry - Howard - Lambek - Scott correspondence also corresponds to the logical `exportation' rule.

This paper is based around currying the elementary arithmetic operations of addition and multiplication. Doing so gives families of injective functions indexed by either the natural numbers $\mathbb N$, or the positive naturals $\mathbb N^+$. We study these indexed families of injections together with their generalised inverses, and describe the resulting inverse submonoids of $\IN$, the symmetric inverse monoid on $\mathbb N$, that they generate.

Some of these are familiar and well-studied (e.g. the bicyclic monoid, the polycyclic monoids, etc.). Others are novel; we recover an inverse monoid that describes `monotone partial injections between congruence classes'. This has a minimal generating set (the generators of all prime-order polycyclic monoids), normal forms (based on mixed-radix number systems), and formul\ae\ for composition based on the Chinese remainder theorem. 

We also compare \& contrast these normal forms with Nivat \& Perot's normal forms for polycyclic monoids. This leads to a characterisation of $p$-adic norms and distances in terms of prime-order polycyclic monoids. %, and similarly for the arithmetic inverse monoid.
%We recover several well-known and well-studied inverse monoids, as well as a single inverse monoid that generalises all of these in a very natural way. As part of this, we recover many well-known notions from number theory in an inverse semigroup theoretic guise.  %Our final destination is inverse monoids and groups of congruential functions, in the sense of \cite{JC}. 

%\begin{remark} [A categorical view]%comment for category-theorists] Hard-core category theorists may object to the description above, on the grounds that 
%	Although this paper is phrased algebraically, rather than categorically, the general setting may also be of interest. 
%	In the category $({\bf Set},\times)$ of sets, functions and Cartesian product, Currying arises as the `tensor-hom adjunction' ${\bf Set}(A\times B,C) \ \cong {\bf Set}(A,[B\rightarrow C])$, where $[B\rightarrow C]\in Ob({\bf Set})$ is the hom-set of functions from $B$ to $C$. Our work may be seen as first restricting this to the setting where hom-sets are {\em injective functions}, applying the faithful functor from $\bf Set$ to the category $\bf pInj$ of sets and partial injective functions, then considering closure under composition and generalised inverses.
%	\end{remark}

\section{Elementary definitions}
We briefly reprise some basics of inverse semigroup theory, partly in order to fix notation and terminology. Readers familiar with inverse semigroups are invited to note our slightly non-standard notation for generalised inverses, identities, and zeros, then to skip forwards to Section \ref{curry-sect}.%, and to refer back to this section for any unfamiliar notation of terminology. 

We assume familiarity with the definitions of monoids and semigroups. The following is gven in order to fix notation : 

%Although not themselves inverse, free monoids often form the basis for significant inverse monoids.

\begin{definition} The free monoid $X^*$ on a set $X$ is the set of all finite strings of members of $X$ (including the empty string $\epsilon\in X^*$), with composition given by concatenation. The {\bf length homomorphism} $len:X^*\rightarrow (\mathbb N , +)$ is the monoid homomorphism defined by $len(x)=1$, for all $x\in X$. The free monoid construction is functorial, so given a function $\phi:X{\rightarrow} Y$, there exists a corresponding homomorphism of free monoids $\phi^*:X^*{\rightarrow} Y^*$. %Categorically, freeness and the universal property of free monoids is exhibited by an adjunction between the `free monoid' functor, and the `forgetful' functor that maps monoids \& homomorphisms to their underlying sets \& functions.
As a corollary, any free monoid $X^*$ on a finite set $|X|=n$ is isomorphic to $\{ 0,\ldots, n-1\}^*$.  We %somewhat sloppily 
refer to this as {\bf the $n^{th}$ free monoid}, and denote it by $F_n \stackrel{def.}{=} \{ 0,\ldots,n-1\}^*$.
	\end{definition} 

The following definitions and results on inverse monoids may be found in, for example, \cite{MVL}.
\begin{definition}
A monoid $M$ is {\bf inverse} when every element $a\in M$ has a unique {\bf generalised inverse} $a^\ddagger\in M$ satisfing $aa^\ddagger a=a$ and $a^\ddagger a a^\ddagger=a^\ddagger$.  The {\bf symmetric inverse monoid} on a set $X$, denoted $\mathcal I (X)$, is the monoid whose elements are {\bf partial injections} on $X$, i.e. partial functions $f:X\rightarrow X$ satisfying $a(x)=a(y) \ \Rightarrow \ x=y$ when $a(x)$ and $a(y)$ both exist.  Composition is the usual composition of partial functions, and generalised inverses are given by $a^\ddagger(y)=x$ iff $a(x)=y$.
%The {\bf diagonal representation} of a partial injection $a\in I(X)$ is the subset $\Delta_a =\{ (y,x) : y=a(x) \}_{y,x\in X} \subseteq X\times X$.

A {\bf zero element} $\epsilon\in M$ is an absorbing element, so $\epsilon a=\epsilon=a\epsilon$ for all $a\in M$. To avoid a clash of notation with the arithmetic we study, we avoid the more common notation $0\in M$ for the zero element. Similarly, we denote the {\bf identity element} by $Id\in M$ rather than the potentially misleading $1\in M$.%confusion 
\end{definition}

For any inverse monoid $M$, and set $X$, the following are standard :
%The following basic facts are standard:
\begin{proposition}\label{basics-prop} \ \\
	\begin{enumerate}
		\item The generalised inverse is a self-inverse anti-isomorphism on $M$, so $Id^\ddagger=Id$,  $(ab)^\ddagger = b^\ddagger a^\ddagger$, and $\left( a^\ddagger\right)^\ddagger=a$ for all $a,b\in M$.
		\item $M$ is isomorphic to a monoid of partial injections (This is the Wagner-Preston representation theorem \cite{WP} -- the inverse semigroup theoretic analogue of Cayley's theorem).
		\item For arbitrary $a\in M$, the elements $a^\ddagger a$  and $aa^\ddagger$ are idempotent (the {\bf initial} and {\bf final idempotents} of $a$ respectively), and all idempotents are of this form.
		\item Idempotents are self-inverse, and all idempotents commute; $e^\ddagger=e$ and $ef=fe$, for all $e^2=e,f^2=f\in M$.
		\item Idempotents may be ``passed through'' %arbitrary 
		elements; given $e^2=e,a\in M$ then $f=aea^\ddagger $ is idempotent, and satisfies $ae=fa$.
		\item The idempotents of $M$ form a meet semilattice, denoted $E(M)$, with meet given by composition and the identity as top element.
		\item The idempotents of $I(X)$ are precisely the partial identities on $X$, that is $Id_S$ for some $S\subseteq X$ given by $Id_S(x)=\left\{\begin{array}{lr} x & x\in S \\ \bot & \mbox{otherwise.}\end{array}\right.$  Idempotents compose as $Id_SId_T=Id_{S\cap T}$ for all $S,T\subseteq X$, and the initial and final idempotents of a partial injection are the partial identities on its domain / image respectively.
		\item Every symmetric inverse monoid has a zero element $\epsilon$, given by the nowhere-defined partial function $\epsilon=Id_{\emptyset }$.
		\end{enumerate}
	\end{proposition}

\section{Currying elementary arithmetic}\label{curry-sect}
We now introduce the core objects of study : partial injections defined by currying elementary arithmetic operations. 
 
\begin{definition}
	Let us denote the natural numbers by $\mathbb N$, and the non-zero naturals by $\mathbb N^+$.
For all $x\in \mathbb N$, and $y\in \mathbb N^+$ let us denote by $+_x,\times_y\in \mathcal I(\mathbb N)$ the injective functions that arise by currying, so {\bf curried addition} and {\bf curried multiplication} are given by, respectively, $+_x(n) = n+x$ and $\times_y(n) = ny$, for all $n\in \mathbb N$. 
Their generalised inverses are the partial injections 
\[ 
+^\ddagger_x(n) = \left\{ \begin{array}{lr}n-x & n\geq x \\ & \\  \bot & \mbox{otherwise,} \end{array}\right. 
\ \  \ \
\times^\ddagger_y(n) = \left\{ \begin{array}{lr} \frac{n}{y} & n \ (\md  \ y) =0 \\ & \\ \bot & \mbox{otherwise.} \end{array}\right.  \] 
which we describe as {\bf non-negative subtraction}, and {\bf whole-number division}, respectively.

Commutativity of addition \& multiplication, and the  
distributivity of multiplication over addition, give the following identities : 
\begin{enumerate}
	\item $+_a+_b=+_b+_a=+_{+_a(b)}$
	\item $\times_b\times_a=\times_a\times_b = \times_{\times_a(b)} $
\item $\times_y+_x \ = \ +_{\times_y(x)}\times_y$
%\item $
\end{enumerate}
%We use these to define 
\end{definition}

%We may use these to define natural generalisations of the `length' homomorphism on a free monoid.

%\begin{definition}
%	Given an arbitrary non-empty set $X$, the well-known {\bf length homomorphism} $| \ |:X^*\rightarrow (\mathbb N, +)$ is the unique monoid homomorphism defined by $|x|=1$, for all $x\in X$.  Analogously, given $a\in \mathbb N$ and $b\in \mathbb N^+$, we define the {\bf $a$-additive} and {\bf $b$-multiplicative length homomorphisms} $|\ |_{+_a} , |\  |_{\times b}:X^*\rightarrow \IN$ to be the unique homomorphisms that map all generators of $X^*$ to, respectively, $+_a,\times_b\in \IN$.
%	\end{definition}

We demonstrate that the inverse submonoids of $\IN$ generated by $\{+_x\}_{x\in \mathbb N}$ and $\{ \times_x \}_{x\in \mathbb N^+}$ are well-known.% inverse monoids. 

\subsection{The bicyclic monoid as curried addition}
The following definition was first written down by E. Lyapin in \cite{EL}, but had previously been studied by other semigroup theorists including Clifford, Preston, \& Rees (see \cite{CH} for a historical overview).
\begin{definition}
	The {\bf bicyclic monoid} $\mathcal B$ is the inverse monoid with a single generator $s$ subject to a single relation $s^\ddagger s=Id$ (Note this is a one-sided inverse, so $Id\neq ss^\ddagger $).  The elements of $\mathcal B$ may be given a normal form as pairs of natural numbers, with composition given by :  
	\[ (d,c) (b,a) \ = \   \left(d + [b\monus c], [c\monus b]+a \right)  \ \ \ \forall \ (d,c),(b,a)\in \mathbb N \times \mathbb N  \] 
	where the {\bf monus} $\monus$ is defined by $y\monus x = \left\{ \begin{array}{lr} y-x & x\leq y \\ 0, & \mbox{otherwise.}\end{array}\right.$% Elements of $\mathcal B$ have normal forms as pairs of natural nmbers 
	\end{definition}
\noindent The following is well-known (see, for example, \cite{MVL}).
\begin{proposition}
	The inverse submonoid of $I(\mathbb N)$ generated by $+_1$ is isomorphic to $\mathcal B$.  
	\end{proposition}

\begin{corollary}\label{Bembed-corol}
The inverse submonoid of $I(\mathbb N)$ generated by $\{+_x\}_{x\in \mathbb N}$ is isomorphic to the bicyclic monoid $\mathcal B$.
\end{corollary}
\begin{proof}It suffices to note that $+_k=\left(+_1\right)^k$ for all $k>0$, and $+_0=Id$. We may then identify normal forms as $[b,a]=+_{b}+_a^\ddagger$.
	\end{proof}

\subsection{Leech's monoid as curried multiplication}
In \cite{JLe}, J. Leech introduced a `multiplicative analogue of the bicyclic monoid' defined as follows:
\begin{definition}
	Leech's {\bf multiplicative monoid} $\mathcal L$ is the inverse monoid whose underlying set is $\mathbb N^+\times \mathbb N^+$, with composition given in terms of greatest common divisors, as 
	\[ (m,n)(p,q) \ = \ \left( \frac{mp}{\gcd (n,p)},\frac{nq}{\gcd (n,p)} \right) \]
	generalised inverse given by $(m,n)^\ddagger=(n,m)$, and idempotents given by $\{ [p,p] : p>0\}$.
	(Leech denoted this monoid by $P$ -- we have changed notation to avoid confusion with the closely related polycyclic monoids of Definition \ref{poly-def}).
	\end{definition}

\begin{remark}\label{otherComp-rem}
	Using the elementary number-theoretic relationship between greatest common divisor and least common multiple, 
	\[ np\ = \ \gcd (n,p)\lcm  (n,p) \ \ \forall n,p\in \mathbb N^+ \] 
	we may give an equivalent formula for composition, % within $\mathcal L$, in terms of least common multiples, 
	as 
	\[ (m,n)(p,q) \ = \ \left( \frac{m.\lcm  (n,p)}{n},\frac{q.\lcm  (n,p)}{p} \right) \]
	%We will interchangeably use either form.%, prefer this equivalent form, particularly when dealing with composites of idempotents.
	\end{remark}

Leech's multiplicative monoid has a relationship with multiplication analogous to the relationship between the bicyclic monoid and addition. We first fix some notation.

\begin{definition}
	Given $a>b\in \mathbb N$, we denote the {\bf congruence class of $b$ modulo $a$} by $a\mathbb N +b  \stackrel{def.}{=} \{ a n + b\}_{n\in \mathbb N}\subseteq \mathbb N$, and refer to $a,b$ respectively as the {\bf multiplicative} and {\bf additive coefficients}. We simplify notation to $a\mathbb N$ when the additive coefficient is zero, and write $1.\mathbb N$ simply as $\mathbb N$. 
	Trivially, each such congruence class is a well-ordered countably infinite subset of $\mathbb N$. They are all therefore in bijective correspondence.
	\end{definition}

\begin{proposition}\label{Lembedding-prop}
	The map $\phi:\mathcal L \rightarrow \IN$ given by $\phi([m,n])=\times_m \times_n^\ddagger$ is an injective inverse monoid homomorphism whose image is the inverse submonoid of $\IN$ generated by $\{ \times_n\}_{n\in \mathbb N^+}$. 
	\end{proposition}
\begin{proof}
	Let us (temporarily) denote the inverse submonoid of $\IN$ generated by $\{ \times_n\}_{n\in \mathbb N^+}$ by $\mathcal L'$. Then $\phi([m,n])=\times_m \times_n^\ddagger$ is the unique monotone partial injection that maps $n\mathbb N$ to $m\mathbb N$ and is undefined elsewhere.   The idempotents of $\mathcal L'$ are of the form $\times_p\times_p^\ddagger $ -- i.e. partial identities on the subset $p\mathbb N$, giving the composition of idempotents as : %, for all $\times_n\times_n^\ddagger ,\times_p\times_p^\ddagger \in E(\mathcal L ')$,
	\[ \times_n\times_n^\ddagger \times_p\times_p^\ddagger  = \times_{\lcm  (n,p)}\times_{\lcm  (n,p)}^\ddagger  
	\ \ \ \ \forall \ n,p\in \mathbb N^+
	\]
	Appealing to the alternative formula for composition given in Remark \ref{otherComp-rem}% well-known identity $np=\gcd (n,p)\lcm  (n,p)$, we derive  
	\[ \phi([n,n][p,p]) = \phi([\lcm  (n,p),\lcm  (n,p)] = \phi([n,n])\phi([p,p]) \] 
	Thus $\phi$ establishes an isomorphism between the semilattices of idempotents of $\mathcal L$ and $\mathcal L'$.
	
	Finally, by monotonicity, elements of both $\mathcal L'$ and $\mathcal L$ are uniquely determined by their initial and final idempotents, giving $\phi$ as an isomorphism. Thus $\mathcal L \cong \mathcal L'$. 
\end{proof}

\begin{corollary}\label{lembed-corol}\ \\
	\begin{enumerate}
		\item $\mathcal L$ is isomorphic to the inverse monoid of all monotone partial injections between congruence classes of the form $\{ a\mathbb N : a\in \mathbb N^+ \}$.%The map $\psi : \mathcal L\rightarrow \IN$ defined by $\psi([m,n])=\times_m\times_n^\ddagger$ is an injective inverse monoid homomorphism. 
		\item The set $\{ [1,p] : p \in \pr   \}$ is a minimal generating set for the inverse monoid $\mathcal L$.
		\end{enumerate}
	\end{corollary}

\section{The arithmetic inverse monoid}
The above embedding of $\mathcal L$ into $\IN$ relies solely on mapping between congruence classes whose additive coefficient is zero (i.e. generated by curried multiplication). The embedding of $\mathcal B$ into $\IN$ does not map between congruence classes at all, but is simply curried addition. 

We now introduce an inverse monoid that combines both curried addition and curried multiplication to give monotone mappings between arbitrary congruence classes. 

By contrast with $\mathcal B$ and $\mathcal L$, this is not a well-known inverse monoid.  However, it contains and generalises in a natural manner Nivat and Perot's polycyclic monoids (Definition \ref{poly-def} below), and may be thought of as being generated by the prime-order polycyclic monoids (Section \ref{MR-sect}). %It also has a very close connection with both $p$-adic arithmetic, and the {\em congruential functions} of Conway (Sections \ref{padic-sect} and \ref{conway-sect}). %Definition \ref{cf-def} via its orthogonal closure.

\subsection{Polycyclic inverse monoids}
The {\em polycyclic monoids} were introduced by Nivat and Perot as a natural generalisation of the bicyclic monoid \cite{NP} and have been repeatedly re-discovered in numerous different settings (notably the `dynamical algebra' of the logicians \cite{DR,GOI1,GOI2}, with the equivalence given in \cite{PHD}), the Cuntz algebras of $C^*$ algebra theory and theoretical physics \cite{K}, and the `bracketing language' or `stack algebra' of theoretical computer science \& automata theory \cite{JS}).   
\begin{definition}\label{poly-def} For any set $X$ with $|X|>1$, the {\bf polycyclic monoid} ${\mathcal P}_X$  is the inverse monoid generated by the set $X$, subject to the relations 
	\[ xy^\ddagger = \left\{ \begin{array}{lcr} Id & & x=y \\
	& & 	\\
	\epsilon & & x \neq y \end{array} \right.  \ \ \ \ \ \forall x,y\in X \]
\end{definition}

\begin{notation}[Free and polycyclic monoids]
	Given sets $X,Y$, with $|X|=|Y|$, then ${\mathcal P}_X \cong {\mathcal P}_Y$.  We therefore apply the same convention as for finitely-generated free monoids, and define {\bf the $k^{th}$ polycyclic monoid} ${\mathcal P}_k$ to be the polycyclic monoid generated by $\{ 0,\ldots , k-1 \}$.%, and somewhat sloppily refer to this as ``''.  We apply the same conventions to the free monoid on this set, and denote `the $k$-generator free monoid' by $\Fk=\{0,\ldots,k-1\}^*$.
\end{notation}% 
%\end{definition}

The following key property, taken from \cite{NP}, determines much of the structure and theory of polycyclic monoids : 

\begin{proposition}\label{cf-prop} All polycyclic monoids are congruence-free -- any homomorphic image of ${\mathcal P}_k$ is either isomorphic to ${\mathcal P}_k$, or is the trivial monoid $\{ Id \}$, for all $k>1$. 
\end{proposition}

\subsection{Polycyclic monoids via arithmetic operations}

The following arithmetic embeddings of finite polycyclic monoids into $\IN$ are well-established (see, for example, \cite{PHD,MVL}).

\begin{theorem}\label{Pkembed-thm}
	For all $k>1$, the map $\theta_k:{\mathcal P}_k\rightarrow \IN$ defined on generators by $\theta_k(x)=(+_x\times_k)^\ddagger$ is an injective inverse monoid homomorphism.
	\end{theorem}
\begin{proof}
	For all $x,y<k\in \mathbb N$, by construction, $dom(+_x\times_k)^\ddagger=k\mathbb N +x$ whereas $im(+_y\times_k)=k\mathbb N + y$. When $y\neq x$ these have empty intersection, giving  $\theta_k(x)\theta_k(y)^\ddagger = (+_y\times_k)^\ddagger (+_x\times_k)= \epsilon$. 
	
	Conversely, for all $n\in \mathbb N$, 
	\[ \theta_k(x)\theta_k(x)^\ddagger (n)\ = \ (+_x\times_k)^\ddagger (+_x\times_k)(n) \ =\ \frac{(kn+x)-x}{k} \ =\ n \]
	Therefore, 
	$(+_x\times_k)^\ddagger \left(+_x\times_k \right)=Id\in \IN$. Hence,  for all $x,y<k$, 
	\[ \theta_k(x)\theta_k(y)^\ddagger\ = \ \left\{\begin{array}{lcr} Id & \ & x=y, \\ & & \\ \epsilon & & x\neq y, \end{array}\right. \]
	and so $\theta_k:{\mathcal P}_k\rightarrow \IN$ is an inverse monoid homomorphism. 
	Finally, $\theta_{k}(x)^\ddagger \theta_k(x) \ = \ Id_{k\mathbb N + x} \ \neq \ Id_\mathbb N $, for all $x<k$,
so the congruence-freeness property implies injectivity. %implies that it is also injective.%an embedding.
	\end{proof}

\begin{remark} For all $k>1$, the embedding $\theta_k:{\mathcal P}_k\hookrightarrow \IN$ is a {\bf strong embedding} in that the identity of $\IN$ arises as the supremum (w.r.t. the natural partial order) of the non-identity idempotents under this embedding, so $Id_{\mathbb N}= \sup \{ \theta_k(e) : e^2=e\neq Id \}$.%non-identity idempotents of {
	\end{remark}

\begin{remark}[Novelty and historical background]
	Representations of polycyclic monoids -- in particular, the two-generator case -- as partial injections on the natural numbers  are certainly not novel to this paper. Girard's representation of the `dynamical algebra' (i.e. ${\mathcal P}_2$) found in \cite{GOI1,GOI2}, used precisely the functions $n\mapsto 2n$ and $n\mapsto 2n+1$, along with their generalised inverses\footnote{As an additional complication, Girard did not explicitly consider partiality. Rather, he worked in an separable infinite-dimensional Hilbert space, using partial isometries with the special property that they mapped elements of some orthonormal basis either to other orthonormal basis vectors, or to the null vector. This corresponds to a slightly disguised version of $\IN$ where partiality need not be considered; `undefined' results are modelled by the null vector.  This is of course an example of a categorical construction -- M. Barr's faithful $l_2$ functor \cite{MBa} that takes partial injections on sets to partial isometries on Hilbert spaces.  Numerous authors \cite{SA96,PHD,DR,AHS} rapidly realised that the Hilbert space structure was inessential, and re-wrote Girard's work in the significantly simpler setting of inverse semigroups. Ironically, the Hilbert space structure was later re-introduced (by several of the same authors, using Barr's functor) for applications of Girard's work to quantum computation (for example \cite{TCS2}).}.  The most general setting was given in \cite{PHD,MVL}, where it was observed that any Hilbert-hotel style bijection $\Psi: \mathbb N \uplus\mathbb N \cong \mathbb N$ uniquely determines and is uniquely determined by a strong embedding ${\mathcal P}_2\hookrightarrow \IN$.
\end{remark}

	Based on the above arithmetic realisation of finitely generated polycyclic monoids, we now propose a natural generalisation.

\subsection{Une g\'en\'eralisation des mono\"ides polycycliques}
The monoid we now introduce may be characterised %(as in Theorem \ref{} below) 
as the inverse submonoid of, %partial injections that (uniquely) 
\begin{center}{\bf ``monotone partial injections between congruence classes''.}\end{center} %We introduce it concretely, and prov This takes place using curried addition and multiplication, and their generalised inverses.
(This intuition is formalised in Corollary \ref{character-corol} below).

It arises in a natural way from the polycyclic monoids; the embeddings  $\{ {\mathcal P}_k \hookrightarrow \IN \}_{ k\in \mathbb N^+}$ of in Theorem \ref{Pkembed-thm} arise via currying the generating set described below. 
\begin{definition}
	Given arbitrary $a>b\in \mathbb N$, 
	we denote the unique monotone partial injection whose domain is $a\mathbb N + b$ and whose image is the whole of $\mathbb N$ by 
	$R_{a,b}=(+_b\times_a)^\dagger \in \IN$, so 
	\[ R_{a,b}(n) = \left\{ \begin{array}{lcr} \frac{n-b}{a} & \ & n\ (\md   \  a) = b \\ & & \\ \bot & & \mbox{otherwise.}
	\end{array}\right.
	\]
	By construction, $R_{a,b}=\theta_a(b)$, where $\theta_a:{\mathcal P}_a\hookrightarrow \IN$ is as given in Theorem \ref{Pkembed-thm}; the generating set $\{ \theta_a(b) : b<a \in \mathbb N \}$ of the $a^{th}$ polycyclic monoid is ``The result of currying the partial injection $R_{a,b}$''.
	
	The generalised inverse $R_{a,b}^\ddagger$ is the unique (globally defined) monotone injection $R_{a,b}^\ddagger=+_b\times_a$ that maps $\mathbb N$ to the congruence class $a\mathbb N +b$.  Explicitly, $R^\ddagger_{a,b}(n)=an+b$. %Thus, $R^\ddagger_{c,d}R_{a,b}$ is the unique monotone injection whose domain is $a\mathbb N+b$ and whose image is $c\mathbb N +d$.

We then define the {\bf arithmetic inverse monoid} $\mathcal A $ to be the inverse subsemigroup of $\IN$ generated by $\{ R_{a,b} : a>b\in \mathbb N \}$. Note that $R_{1,0} = Id\in \IN$, so $\mathcal A$ is a monoid.
\end{definition}

The following are immediate from direct calculation, and are relevant for the normal forms we establish in Section \ref{Anf-sect}.

\begin{proposition}
Given $a>b,c>d\in \mathbb N$, then 
\[ R_{c,d}^\ddagger R_{a,b}(n) \ = \ \left\{ \begin{array}{lcr} 
c\left( \frac{n-b}{a}\right) + d & \ \  & n \ (mod \ a)=b \\ & & \\ 
\bot & & \mbox{otherwise.}
\end{array}\right. 
\]
When $R_{c,d}^\ddagger R_{a,b}(n)$ is defined, it is natural to write it in terms of a $(2\times 2)$ determinant, so
$R^\ddagger_{c,d}R_{a,b}(n) \ = %\ \frac{cn-cb+ad}{a} \ =
\ \frac{1}{a} \left( cn + \left| \begin{array}{cc} a & b \\ c & d \end{array} \right| \right)$.
\end{proposition}
\begin{proof}
	These are immediate from the definitions, and very elementary arithmetic manipulation.
	\end{proof}
We now establish some embeddings of the inverse monoids described in Section \ref{curry-sect}. %. These are straightforward, and motivated the definition of $\mathcal A$.
\begin{theorem}\label{Aembeddings-thm} The arithmetic monoid $\mathcal A$ contains :
	\begin{enumerate}
		\item  an isomorphic copy of Leech's $\mathcal L$, generated by $\{ R_{a,0}:a\in \mathbb N^+\}$
		\item for all $a>1$, an isomorphic copy of the polycyclic monoid ${\mathcal P}_a$ generated by $\{ R_{a,b}:b<a\in \mathbb N\}$ 
		\item countably infinitely many distinct isomorphic copies of the bicyclic monoid $\mathcal B$.
		\end{enumerate}
	\end{theorem}
\begin{proof}
	These are simple corollaries of previous results.
	\begin{enumerate}
		\item Recall the embedding $L\hookrightarrow \IN$ of Corollary \ref{lembed-corol} given by $[m,n]\mapsto \times_m\times_n^\ddagger$, and observe that $\times_m\times_n^\ddagger=R_{m,0}^\ddagger R_{n,0}\in \mathcal A$. Thus the inverse submonoid of $\mathcal A$ generated by $\{ R_{n,0}:n>0 \}$ is isomorphic to $\mathcal L$.
		\item This is by construction. The embedding $\theta_a:{\mathcal P}_a\hookrightarrow \IN$ satisfies $\theta_a(b)=R_{a,b}$; thus this embedding factors through $\mathcal A$. % and so the following that . 
		\item Some subtlety is needed in order to give a suitable embedding. The set $\mathbb N+k$ is not a congruence class, for all $k>0$, so the embedding $\mathcal B \hookrightarrow \IN$ implicit in Corollary \ref{Bembed-corol} does not factor through $\mathcal A \hookrightarrow \IN$. Instead,  we rely on the `exponential' embeddings, $p^{(\ )} : \mathcal B \rightarrow \mathcal L$ given in Leech's original paper \cite{JLe} as $p^{([b,a])}=[p^b,p^a]$, for arbitrary prime\footnote{The result of \cite{JLe} is given in the case where $p$ is prime, although the proof does not rely on primality. However,  it forms part of a proof of a stronger result on embedding countably infinite products of $\mathcal B$ into $\mathcal L$ where primality -- or at least, a countably infinite set of co-prime naturals -- is necessary.}  $p\in {\tt primes}$. Composing with the embedding of part 1. above gives our result.
	\end{enumerate}	
	\end{proof}

	We now give relations between these generators to characterise elements of $\mathcal A$, and give a minimal generating set, normal forms, and formul\ae\ for composition. These rely on a staple of number theory :

\begin{theorem}[The Chinese Remainder Theorem]
Given $a>b,c>d\in \mathbb N$, then $(a\mathbb N+b)\cap (c\mathbb N + d)$ is either 
\begin{itemize}
	\item of the form $\lcm  (a,c)\mathbb N +y$, for some $y<\lcm  (a,c)$ when \[ |b-d| \ \in \ \gcd (a,c)\mathbb N \]
	\item the empty set, otherwise.% when $(b\monus d) + (d\monus b) \ \notin \ \gcd (a,c)\mathbb N$.
\end{itemize}
\end{theorem}
\begin{proof}
	This is given in every introduction to number theory (e.g. \cite{WS}), and the (entirely constructive) procedure for finding $y<\lcm  (a,c)\in \mathbb N$ when a solution exists is a common component of introductory number theory courses -- at least in the simple case where $a,c\in \mathbb N$ are co-prime, in which case a solution necessarily exists.
\end{proof}

%\end{document}

\begin{proposition}\label{bd-prop}
	For all $d<c$ and $b,b'<a\in \mathbb N$,
	\begin{enumerate}
		\item All elements $f\in\mathcal A$ are monotone, so $f(x)\leq f(y)$ iff $x\leq y$, provided $f(x),f(y)$ are defined.
		\item $R_{a,b}^\ddagger R^\ddagger_{c,d}= R^\ddagger_{ac,ad+b}$, and hence $R_{c,d}R_{a,b}=R_{ac,ad+b}$.
		\item $R^\ddagger_{c,d}R_{a,b}$ is the unique monotone injection that maps $a\mathbb N+b$ to $c\mathbb N +d$, and is undefined elsewhere, and as a special case, $E_{a,b}=R^\ddagger_{a,b}R_{a,b}$ is idempotent.
		\item 
		\[ R_{c,d}R_{a,b}^\ddagger = \left\{ \begin{array}{lcr} R_{w,x}^\ddagger R_{u,v} & & |b-d| \ \in \ \gcd (a,c)\mathbb N \\ & & \\ \epsilon & & \mbox{otherwise,}
		\end{array}\right. \]
		 where, in the non-empty case,
		 $(a\mathbb N +b) \ \cap \ (c\mathbb N + d) \ = \ \lcm  (a,c)\mathbb N +r$, 
		 giving
		 \begin{itemize}
		 	\item $u =\frac{c}{\gcd (a,c)}$ and $v=\frac{r-b}{a}$
		 		\item $w=\frac{a}{\gcd (a,c)}$ and $x=\frac{r-d}{c}$
		 		\end{itemize}
	\end{enumerate}
\end{proposition}
\begin{proof} \ \\
	\begin{enumerate}
		\item We observe all generators to be partial monotone, and note that this property is preserved under both composition and generalised inverses.
		\item By definition, for all $n\in \mathbb N$,
		\[ R^\ddagger_{a,b}R^\ddagger_{c,d}(n) \ = \ R^\ddagger_{a,b}(cn+d) \ = \ acn+ad+b \ = \ R^\ddagger_{ac,ad+b} \]
		\item This is essentially by construction; $R_{a,b}$ is the unique monotone injection mapping $a\mathbb N+b$ to $\mathbb N$, and $R^\ddagger_{c,d}$ is the unique monotone injection mapping $\mathbb N$ to $c\mathbb N + d$. As a special case, $E_{x,y}=R^\ddagger_{x,y}R_{x,y} = Id_{x\mathbb N + y}$.
		\item This is an appeal to the Chinese Remainder Theorem. 
		The relevant domains and images are $im(R^\ddagger_{a,b})=a\mathbb N +b$, and $dom(R_{c,d})=c\mathbb N + d$. Their intersection is non-empty precisely when $|b-d|\in \gcd (a,c)\mathbb N$, in which case
		\[ R_{c,d}R_{a,b}^\ddagger \ = \ R_{c,d}E_{\lcm  (a,c),q} R_{a,b}^\ddagger \ = \ R_{c,d}E_{\lcm  (a,c),q} E_{\lcm  (a,c),q} R_{a,b}^\ddagger \]
		Observing that the domain of $R_{c,d}E_{\lcm  (a,c),q}$ is identical to the image of $E_{\lcm  (a,c),q} R_{a,b}^\ddagger$, we may pass these idempotents through the elements (Part 5. of Proposition \ref{basics-prop}), and derive
	\[ dom\left(E_{\lcm  (a,c),r} R_{a,b}^\ddagger\right) \ = \ \frac{\lcm   (a,c)}{c} \mathbb N + \frac{r-b}{a} \]
and similarly
 	\[ im \left( R_{c,d}E_{\lcm  (a,c),r} \right) \ = \ \frac{\lcm  (a,c)}{c} \mathbb N + \frac{r-d}{c} \]
 		Appealing to the identity $ac=\lcm  (a,c)\gcd (a,c)$ gives that $R_{c,d}R_{a,b}^\ddagger$ is the unique monotone injection whose domain is $\frac{c}{\gcd (a,c)}\mathbb N + \frac{r-b}{a}$ and whose image is $\frac{a}{\gcd (a,c)}\mathbb N + \frac{r-d}{c}$. However, by part 3. above, this is necessarily the partial injection $R_{w,x}^\ddagger R_{u,v}$, as required.
		\end{enumerate}
	\end{proof}

The following characterisation of $\mathcal A$, along with normal forms for elements (Section \ref{Anf-sect}), follows directly from the above results.

\subsection{Normal forms for $\mathcal A$}\label{Anf-sect}
As a corollary of Proposition \ref{bd-prop}, $\{ R_{c,d}^\ddagger R_{a,b} : a>b,c>d\in \mathbb N \} \cup \{ \epsilon \}$ is in fact closed under composition -- explicitly, giving normal forms for elements of $\mathcal A$. 
We prove this, and describe composition of elements written in normal form.
%\pagebreak 
\begin{theorem}\label{Anf-thm}\ \\ \begin{enumerate}
		\item 
	Every non-zero element of $\mathcal A$ may be written as $R^\ddagger_{c,d}R_{a,b}$ for some $c>d,a>b\in \mathbb N$. 
	\item Composition of normal forms is given by 
	$\left( R^\ddagger _{g,h} R_{e,f} \right) \left( R^\ddagger_{c,d}R_{a,b} \right)  \ = \ $
	\[ \left\{ \begin{array}{lcr}
		R^\ddagger_{gw,gx+h}R_{au,av+b} & \ \ & |d-f|\in \gcd (c,e)\mathbb N \\ & & \\ \epsilon & & \mbox{otherwise,}
	\end{array} \right.
		\]
		where, in the non-empty case, $(c\mathbb N + d) \cap (e\mathbb N + f ) = \lcm  (c,e)\mathbb N + r$,  
		giving 
		\[ \begin{array}{ccc} u = \frac{e}{\gcd (c,e)} & & v=\frac{r-d}{c} \\ & & \\  
								w=\frac{c}{\gcd (c,e)} & & x=\frac{r-f}{e} 
									\end{array}
									\]
		Explicitly, %substituting in the values of $u,v,w,x\in \mathbb N$  gives  
		\[ \left( R^\ddagger _{g,h} R_{e,f} \right) \left( R^\ddagger_{c,d}R_{a,b} \right) \ = \ 
		R^\ddagger_{\frac{gc}{\gcd (c,e)},\frac{g(r-f)}{e}+h} R_{\frac{ae}{\gcd (c,e)},\frac{a(r-d)}{c}+b} \]
		Using the relationship $mn=gcd(m,n)lcm(m,n)$, we may also write this as 
		\[ \left( R^\ddagger _{g,h} R_{e,f} \right) \left( R^\ddagger_{c,d}R_{a,b} \right) \ = \ 
		R^\ddagger_{\frac{g .\lcm (c,e)}{e},\frac{g(r-f)}{e}+h} R_{\frac{a.\lcm(c,e)}{c},\frac{a(r-d)}{c}+b} \]
		\end{enumerate}
			\end{theorem}
		\begin{proof} \ \\
			\begin{enumerate}
				\item 
			To see that elements may indeed be written in normal form, consider an arbitrary string of generators / generalised inverses thereof. Part 2. of Proposition \ref{bd-prop} demonstrates that every pair of generators $R_{\_,\_}R_{\_,\_}$ may be replaced by a single generator $R_{\_,\_}$, and similarly for generalised inverses of generators;  every $R_{\_,\_}^\ddagger R_{\_,\_}^\ddagger $ may be replaced by a single $R_{\_,\_}^\ddagger$. Now consider an adjacent pair of the form $R_{\_,\_} R_{\_,\_}^\ddagger$. Part 3. of Proposition \ref{bd-prop} allows us to replace this either by the zero element $\epsilon$, or by a pair of the form $R^\ddagger_{\_,\_}R_{\_,\_}$.
			
			Thus we have a confluent, strictly reducing, rewriting schema on strings of generators that will arrive at either the zero arrow, or a (unique) word of the form $R^\ddagger_{\_,\_}R_{\_,\_}$.
			\item 
			To give an explicit formula for composites in normal form, consider $R^\ddagger _{g,h} R_{e,f} , R^\ddagger_{c,d}R_{a,b} \in \mathcal A$. Let us assume that this composite is non-zero : equivalently, $|d-f|\in \gcd (c,e)\mathbb N$, so there exists some $r\in \mathbb N$ such that 
			\[ \left( e\mathbb N + f \right) \cap  (c\mathbb N +d ) \ = \ \lcm  (c,e)\mathbb N + r \]
			Part 3. of Proposition \ref{bd-prop} gives $R_{e,f}R^\ddagger_{c,d} = R^\ddagger_{w,x}R_{u,v}$ where 
				\[ u = \frac{e}{\gcd (c,e)} \ \ , \ \  v=\frac{r-d}{c} \ \ , \ \  
			w=\frac{c}{\gcd (c,e)} \ \ , \ \  x=\frac{r-f}{e}
			\]
			Two applications of Part 2. of Proposition \ref{bd-prop} then give
			\[ \left( R^\ddagger _{g,h} R_{e,f} \right) \left( R^\ddagger_{c,d}R_{a,b} \right) \ = \ R^\ddagger_{gw,gx+h}R_{au,av+b} \]
			as required.
		\end{enumerate}
	\end{proof}

		\begin{corollary}\label{character-corol}
		%	\begin{lemma}
			%	\item Every non-zero element $f\in A$ satisfies $f =  R^\ddagger_{d,c}R_{a,b}$ for some $c>d,a>b\in\mathbb N$.
			%	\item $R^\ddagger_{d,c}R_{a,b}=R^\ddagger_{d',c'}R_{a',b'}$ iff $a=a',b=b',c=c',d=d'$
				%\item 
				The arithmetic inverse monoid $\mathcal A$ may be characterised as the inverse monoid whose non-zero elements map between congruence classes in a monotone manner. Consequently, the non-zero idempotents of $\mathcal A$ are precisely partial identities on congruence classes. 
	%		\end{corollary}
		\end{corollary}
	%	\begin{proof}
			This follows from Theorem \ref{Anf-thm} above, and Parts 1. and 3. of Proposition \ref{bd-prop}.
	%		Part 1. follows directly from Parts 2. and 4. of Proposition \ref{bd-prop}.  Part 2. follows from part 3. of Proposition \%ref{bd-prop}, as does Part 3. Part 4. simply relies on Part 3. to establish that these are the only idempotents of $\mathcal A$.
	%	\end{proof}
	
\subsection{Primes, generators, and mixed-radix systems}\label{MR-sect}
We now move on to viewing the arithmetic inverse monoid -- and hence, all monotone partial injections between congruence classes -- as :  \begin{center}{\bf ``generated by all prime-order polycyclic monoids''}.\end{center} The key to this is the following lemma, which demonstrates a close connection between composition in $\mathcal A$, and mixed-radix counting systems.

\begin{lemma}\label{mixedradix-lem}
	Consider some indexed set of pairs of natural numbers, $\{ a_j > b_j\in \mathbb N \}_{j=0\ldots n}$. Then $
	R_{a_n,b_n}R_{a_{n-1},b_{n-1}} \ldots R_{a_0,b_0} \ = \ R_{A , B }
	$
	where %$A=a_na_{n-1}\ldots a_0$, and $B= b_0+a_0b_1 +a_0a_1 b_2 + \ldots  + a_0a_1\ldots a_{n-1}b_n$. 
%	More concisely, 
	\[ A=\prod_{j=0}^n a_j \ \ \ \mbox{ and } \ \ \  B=b_0 + \sum_{j=1}^n\left( b_j \prod_{i=0}^{j-1}a_i \right) \] 
	%provided we define the product of the empty indexed set to be $1\in \mathbb N$.
	\end{lemma}
\begin{proof}
	This follows by induction; the first step is given by Part 2. of Proposition \ref{bd-prop}, and the induction step follows by direct calculation.
	\end{proof}

\begin{remark}[Composition in $\mathcal A$ via mixed-radix counting systems] \label{mixedradix-rem}
	The link between composition within the arithmetic inverse monoid, and mixed-radix counting is then almost immediate. 
	
	Let $B\in \mathbb N$ be given by treating a string $b_nb_{n-1}b_{n-2}\ldots b_{0}$ as the representation of a natural number in a mixed-radix counting system, where the respective columns are labelled by $a_n,a_{n-1},a_{n-2},\ldots ,a_0$, giving
	 	\[ B \ = \ \begin{array}{|c|c|c|c|c|}
	 \hline
	 \mbox{Base } \ a_n & \mbox{Base } \ a_{n-1} & \ldots & \mbox{Base } \ a_1 & \mbox{Base } \ a_0 \\ 
	 \hline 
	 b_n & b_{n-1}  & \ldots & b_1 & b_0   \\ \hline 
	 \end{array} 
	 \]
	 The well-established formula
	 \[ B=b_0 + \sum_{j=1}^n\left( b_j \prod_{i=0}^{j-1}a_i \right) \]
	 (commonly attributed to G. Cantor \cite{GC}) then makes the interpretation as mixed-radix counting immediate. 
	   
	\end{remark}

This gives a direct route into formalising the intuition that $\mathcal A$ is, `generated by all prime-order polycyclic monoids'. The following is a relatively straightforward corollary of Lemma \ref{mixedradix-lem} above. 

\begin{theorem}
	Recall the embedding $\theta_k : {\mathcal P}_k\hookrightarrow \IN$ of Theorem \ref{Pkembed-thm}. The arithmetic monoid $\mathcal A$ is the inverse submonoid of $\IN$ generated by : 
	\[ {\mathcal GP}   \ = \{ \theta_p(a) : a<p \}_{p\in \pr  } \ \subseteq \ \IN \]
	(i.e. the image of all the generators of the $p^{th}$ polycyclic monoid, for all prime $p$),
	and no proper subset of this set generates $\mathcal A$.
\end{theorem}
\begin{proof}
	By construction, $\theta_{p}(a) = R_{p,a}$ so the set ${\mathcal GP}  $ is a subset of the generating set of $\mathcal A$. Now consider some generator $R_{x,y}$ of $\mathcal A$; we demonstrate that it may be given as a composite of elements of ${\mathcal GP}  $.  If $x$ is prime, there is nothing to prove. Instead, let us assume that $x={\mathcal P}_1\ldots {\mathcal P}_a$ for some indexed family of primes $\{ {\mathcal P}_j \in {\tt primes} \}_{j=1\ldots a}$ (which is, of course, {\em unique up to re-indexing} or some permutation of $\{ 1,\ldots , j\}$).
	
	As a triviality (based on Part 1. of Theorem \ref{Aembeddings-thm}), when $y=0$, we have $R_{x,0}=R_{{\mathcal P}_1,0}{\mathcal P}_{{\mathcal P}_2,0},\ldots R_{{\mathcal P}_a,0}$. and any permutation of this composite will similarly give $R_{x,0}\in \mathcal A$.  When $y\neq 0$, the ordering assigned to the prime factorisation of $x$ becomes significant. Let us write the natural number $y$ in a mixed-radix counting system as 
	\[ y \ = \ \begin{array}{|c|c|c|c|}
	\hline
	\mbox{Base } \ {\mathcal P}_1 & \mbox{Base } \ {\mathcal P}_2 & \ldots & \mbox{Base } \ {\mathcal P}_a \\ 
	\hline 
	q_1 & 		q_2 		& \ldots & q_a  \\ \hline 
	\end{array} 
	\]
	where $0\leq q_j <{\mathcal P}_j$, for all $j=1\ldots a$. We observe that every natural number in $\{ 0,\ldots, x-1\}$ can be written in this form.
	
	Direct calculation, based on Part 2. of Proposition \ref{bd-prop} then demonstrates that 
	$R_{x,y}\ =\ R_{{\mathcal P}_1,q_1}R_{{\mathcal P}_2,q_2}\ldots R_{{\mathcal P}_a,q_a}$ 
	giving the required factorisation into elements of ${\mathcal GP}  $.
	
	Finally, to demonstrate that ${\mathcal GP}  $ is a minimal generating set, it suffices to consider idempotents. For $p\in {\tt primes}$, and $0\leq x<p\in \mathbb N$, the initial idempotent of the corresponding generator satisfies $R^\ddagger_{p,x}R_{p,x}=Id_{p\mathbb N + x}$. We observe that no $p\in {\tt primes}$ may be written as the least common multiple of any set other than $\{1, p\}\subseteq \mathbb N$. 
\end{proof}

\subsection{Nivat \& Perot's normal forms, arithmetically}
We have given normal forms, and formul\ae\ for their composition,  for $\mathcal A$, together with an interpretation as the inverse monoid generated by the prime-order polycyclic monoids. We now compare these normal forms with the usual normal forms for polycyclic monoids with the normal forms for polycyclic monoids.  %our number-theoretic formul\ae\ for composition. 

In \cite{NP}, Nivat and Perot gave normal forms for elements of polycyclic monoids as pairs of words in the free monoid over the generating set. Their composition was based on `cancellation of matching substrings'.
The following is taken from \cite{NP} : 

\begin{theorem}\label{NFcomp-thm} Consider the polycyclic monoid ${\mathcal P}_k$, for some $k>1$, together with words 
	$x^\ddagger w , v^\ddagger u\in {\mathcal P}_k$ where $x,w,v,u\in F_k$. 
	The composite $x^\ddagger wv^\ddagger u\in {\mathcal P}_k$ is given by one of the following three possibilities :
	\begin{enumerate}
		\item $w$ is of the form $w=rv$, for some $r\in \{ {\mathcal P}_0,\ldots,{\mathcal P}_{k-1}\}^*$, in which case
		\[ x^\ddagger wv^\ddagger u\ = \ x^\ddagger rvv^\ddagger u = x^\ddagger(ru) \]
		\item $v$ is of the form $v=sw$, for some $s\in F_k$, in which case 
		\[ x^\ddagger wv^\ddagger u = x^\ddagger w (sw)^\ddagger u = x^\ddagger ww^\ddagger s^\ddagger u = x^\dagger s^\ddagger u = (sx)^\ddagger u \]
		\item Neither 1. nor 2. hold, in which case $x^\ddagger wv^\ddagger u = \epsilon$.
	\end{enumerate}
\end{theorem}
\begin{proof}
	This was first proved in \cite{NP}. Parts 1. and 2. are, as shown above, almost immediate. Part 3. follows by an induction argument demonstrating that, when neither 1. nor 2. hold, the composite $wv^\ddagger$ will necessarily contain some $ij^\ddagger$ with $i\neq j$, giving  $ij^\ddagger=\epsilon$.
\end{proof}

\begin{corollary} All elements of ${\mathcal P}_k$ may be written in the form $v^\ddagger u$, for some $u,v\in \mathcal F_k$.
	\end{corollary}

\begin{corollary}\label{NPNF-corol} For $k\geq 2$, elements of ${\mathcal P}_k$ have normal forms as members of $\mathcal F_k\times \mathcal F_k$, with composition given by, for all $(x,w),(v,u)\in \mathcal F_k\times \mathcal F_k$,
	\[ (x,w)(v,u) \ = \ \left\{ \begin{array}{lcr}
								(x,ru) & \ & w=rv, \\ & & \\
								(sx,u) & & v=sw, \\  & & \\
								(\epsilon,\epsilon) & & \mbox{otherwise.} \end{array}\right. \]
								\end{corollary}

We interpret such normal forms arithmetically, via the embeddings $\theta_k:{\mathcal P}_k\hookrightarrow \mathcal A$, and compare with the normal forms for elements of $\mathcal A$. As a special case of the close connection between the arithmetic inverse monoid and mixed-radix counting systems described in Lemma \ref{mixedradix-lem} and Remark \ref{mixedradix-rem}, we observe a close connection with representations of natural numbers within positional number systems. % (i.e. the special case of mixed radix counting systems where each column has the same label). 
The following simple corollary / special case of Lemma \ref{mixedradix-lem} will be key :

\begin{corollary}\label{radix-corol}
	Given an indexed family natural numbers $\{ x_i<k \}_{i=0\ldots a-1}$, then  
	\[ R^\ddagger_{k,x_0} R^\ddagger_{k,x_1} \ldots R^\ddagger_{k,x_{a-1}} \ = \ R^\ddagger_{k^a,N} \]
	where $N=\Sigma_{i+j=a-1} k^i x_j$. %Similarly, by considering generalised inverses,
%	\[ \theta^{1,0}_{k,x_0} \theta^{1,0}_{k,x_1} \ldots \theta^{1,0}_{k,x_{a-1}} \ = \ \theta_{k^a,N}^{1,0} \]
\end{corollary}
%\begin{proof} This is precisely a special case of Lemma \ref{mixedradix-lem}.
%	\end{proof}

\begin{remark}\label{basis-rem}
	Our overall claim is that Nivat and Perot's normal forms correspond to ``the special case of mixed-radix counting, where the same base labels each column'' --- e.g. our familiar decimal system for ${\mathcal P}_{10}$, and binary or hexadecimal for ${\mathcal P}_2$ or ${\mathcal P}_{16}$ respectively. The above formula simply gives $N$ as the string $x_0x_1\ldots x_{a-1}$, considered as a number written out in base $k$ (possibly with leading zeros). We use this to give alternative equivalent normal forms for elements of polycyclic monoids.%\footnote{or indeed, normal forms for $k$-adic numbers, when $k$ is prime.}; it is simply $x_{a-1}x_{a-1}\ldots x_0$, considered as a number written in base $k$. 
	%This gives a neat way of describing the embedding $\psi_k: {\mathcal P}_k\hookrightarrow \mathcal A $, in terms of normal forms.%via normal forms.% of elements of (finitely generated) polycyclic monoids.
\end{remark}

%We now expand on Remark \ref{basis-rem}, and formalise the observation that embeddings of normal forms arise from considering words in the free monoid ${\mathcal P}_k$ as numbers written in base $k$.

\begin{definition} Given some $k\in \mathbb N^+$, there is an obvious way of interpreting strings over $\{ 0,\ldots,k-1\}$  (i.e. elements of the free monoid $\mathcal F_k$ as natural numbers written in base $k$ (again, possibly with leading zeros).  We define the {\bf numeric interpretation} to be the function $ num:\Fk \rightarrow \mathbb N$ given by 
\begin{itemize}
	\item $num(\epsilon)=0\in \mathbb N$,
	\item $num(x)=x\in \mathbb N$, 
	\item $num(wx)=\times_k(num(w))+num(x)$,
	\end{itemize}
for all $x\in \{ 0,\ldots ,k-1\}$ and $w\in \{ 0,\ldots ,k-1\}^*$.
Although the intuition is straightforward, this is not well-behaved in a semigroup-theoretic sense. It is simply a surjective, but not injective, function, and there is no reasonable sense in which it is a monoid homomorphism.
\end{definition}

\begin{theorem} Let us fix some arbitrary $k>1$, and consider some element in normal form $v^\ddagger u\in {\mathcal P}_k$, so $u,v\in \Fk$. Then 
	\[ \theta_k (v^\ddagger u) = R^\ddagger_{k^{len(v)},num(v)}R_{k^{len(u)},num(u)} \] 
\end{theorem}
\begin{proof} This follows directly from Lemma \ref{radix-corol} above, and the fact that $\theta_k:{\mathcal P}_k\rightarrow \mathcal A $ is an injective inverse monoid homomorphism. % When $k$ is prime, we also observe the apparent structural connection with normal forms of $k$-adic numbers.
\end{proof}

\begin{corollary} The image of $\theta_k:{\mathcal P}_k\hookrightarrow \mathcal A $ is the set 
	\[ \theta_k\left( {\mathcal P}_k\right) \ = \ \left\{ R^\ddagger_{{k^y},Y}R_{k^x,X} : Y<k^y,X<k^x \right\} \ \subseteq \ \IN \]
\end{corollary}
\begin{proof}
	This follows directly from Corollary \ref{radix-corol} above, and Nivat \& Perot's normal forms (Corollary \ref{NPNF-corol}).
\end{proof}

Based on this, we may give a number-theoretic reinterpretation of  Nivat \& Perot's normal forms and formul\ae\ for composition. We first establish some preliminary definitions and notation.

\begin{definition}
	For all $k>1$, we define the monoid of {\bf $\bf k^{(\ )}$-bounded naturals} $\KBN$ to have underlying set $ \{ (m,n) : n<k^m \in \mathbb N \} \ \cup \ \{ \epsilon\}$ and composition 
	\[ (d,c)\cdot (b,a) = (d+c,k^bc+a) \ \ \mbox{ and } \ \ \epsilon g = g = g\epsilon \ \forall g\in \KBN \]
	Direct calculation demonstrates that this is a $k$-generator free monoid. 
	\end{definition}

\begin{lemma} $\KBN \cong \mathcal F_k$.
	\end{lemma}
\begin{proof} Define $\mu_k: \Fk \rightarrow \KBN$ by $\mu_k(w) = (len(w),num(w))$; this is a monoid isomorphism.
	\end{proof} 
The crucial `matching substrings' condition from Nivaat \& Perot's composition of normal forms for polycyclic monoids (Corollary \ref{NPNF-corol}) interprets arithmetically within $\KBN$
\begin{definition}
	Given $(b,a),(y,x)\in \KBN$, call $(y,x)$ a {\bf $\bf k$-residue} of $(b,a)$ when $y\leq b$ and $a\ (\md   \ k^y)=x$. Then define the {\bf cancellation} to be a partial function $\_ \backslash \_ : \KBN\times\KBN\rightarrow \KBN$ where $(b,a)\backslash (y,x)$ is defined iff $(y,x)$ is a $k$-residue of $(b,a)$, in which case 

	\[ (b,a) \backslash  (y,x) \ = \ 
	\left\{ \begin{array}{lcr} \left( b-y , \frac{a-x}{k^y} \right) & \ & b>y \\
									& & \\
									\epsilon & & b=y \end{array} \right. \]
	\end{definition}

	\begin{proposition}\label{modularcancel-prop}
		Given $w,v\neq \epsilon \in \Fk$, then $\mu_k(v)\in \KBN$ is a $k$-residue of $\mu_k(w)\in \KBN$ iff there exists some $r\in \Fk$ such that $w=rv$, in which case $\mu_k(r)= \mu_k(w) \backslash \mu_k(v)$.
		\end{proposition} 
	\begin{proof}
		 \ \\
		 $(\Leftarrow)$ Consider $v,r\in \Fk$. As $\mu_k$ is homomorphic,% isomorphism
		 \[ \mu_k(rv) \ = \ \mu_k(r)\mu_k(v) \ =\  (len(r),num(r))(len(v),num(v)) \]
		 By definition of composition in $\KBN$,
		 \[ \mu_k(rv) \ = \ \left( len(r)+len(v), k^{len(v)} num(r)+num(v) \right) \]
		 Trivially, $(k^{len(v)}num(r)+v) \ \left(\md  \ k^{len(v)}\right) = v$, so $\mu_k(v)$ is a $k$-residue of $\mu_k(w)$. Finally, observe that 
		 \[ \begin{array}{rcl}
		  \mu_k(w)\backslash \mu_k(v) & = & \left( \mu_k(r)\cdot \mu_k(v) \right) \backslash \mu_k(v) \\
		  && \\
		  	& = &  \left( len(r)+len(v), k^{len(v)} num(r)+num(v) \right) \backslash \mu_k(v) \\
		  	&& \\
		  	& = & \left(len(r)+len(v)-len(v) , \frac{k^{len(v)}num(r)+num(v) - num(v)}{k^{len(v)}}\right)  \\
		  	&& \\
		  	& = & (len(r),num(r) ) \end{array}
		  	\]
		  	as required.
		 	\\
		 	$(\Rightarrow)$ This follows directly, since $\mu_k$ is an isomorphism. \end{proof}

\begin{corollary}
	The set $\KBN\times\KBN$, equipped with the composition\\  $\left[ (x,X),(w,W) \right] \star \left[ (v,V),(u,U) \right] \ =  \ $
	\[
	\left\{ \begin{array}{lcr}
	\left[ (x,X) , \left( (w,W) \backslash (v,V) \right) \cdot (u,U) \right] & \ \ & (v,V) \mbox{ is a $k$-residue of } (w,W) \\ 
		& & \\
	\left[ \left( (v,W) \backslash (w,W) \right) \cdot(x,X) ,  (u,U) \right] & \ \ & (w,W) \mbox{ is a $k$-residue of } (v,V) \\ 
		&	& \\
	\left[ \epsilon,\epsilon \right] & & \mbox{otherwise.}
	\end{array} \right. \]
	is isomorphic to ${\mathcal P}_k$.
	\end{corollary}
\begin{proof}
	This follows from Nivat \& Perot's normal forms for polycyclic monoids, and the arithmetic interpretation of string-cancellation given in Proposition \ref{modularcancel-prop} above.
	\end{proof}

The translation into elements of $\mathcal A$ is similarly straightforward.

\begin{corollary}
	Given elements of $\mathcal A$ in normal form, $R^\ddagger_{k^x,X}R_{k^w,W}$ and $R_{k^v,V}^\ddagger R_{k^u,U}$, then either 
	\begin{enumerate}
		\item $(v,V)$ is a $k$-residue of $(w,W)$, in which case 
		\[ \left(R^\ddagger_{k^x,X}R_{k^w,W}\right)\left(R_{k^v,V}^\ddagger R_{k^u,U}\right) \ = \ R_{k^x,X}^\ddagger R_{k^s,S} , \]
		where $(s,S)= \left( (w,W) \backslash (v,V) \right) \cdot (u,U)$
		\item $(w,W)$ is a $k$-residue of $(v,V)$, in which case 
		\[ \left(R^\ddagger_{k^x,X}R_{k^w,W}\right)\left(R_{k^v,V}^\ddagger R_{k^u,U}\right) \ = \ R_{k^t,T}^\ddagger R_{k^u,U} , \]
		where $(t,T)= \left( (v,W) \backslash (w,W) \right) \cdot(x,X) ,  (u,U) $
		\item $\epsilon$ when neither 1. nor 2. hold. 
	\end{enumerate}
	\end{corollary}

\subsection{From polycyclics to $p$-adics}\label{padic-sect}
	There appears to be a close connection between the minimal generating set of $\mathcal A$ (i.e. all prime-order polycyclic monoids) and $p$-adic arithmetic. Comparing the well-established use of polycyclic monoids in studying fractal structures \& self-similarity (e.g. \cite{PHD}) with the well-known fractal structure of $p$-adic arithmetic {e.g. \cite{AMR}), this is perhaps unsurprising. 
	
	We give a brief overview of some very basic definitions. %, and their interpretation as inverse semigroup theory.
	The following may be found in many number theory texts, such as \cite{AMR}.
	\begin{definition}
		Given $p\in {\tt primes}$, the {\bf $p$-order} $ord_p:\mathbb N\rightarrow \mathbb N$ is defined by 
		\[ ord_p(n) \ = \ max \left \{k : n \ (\md \ p^k) = 0 \right\} _{k\in \mathbb N} \]
		i.e. it is the number of occurrences of $p$ in the (unique) prime factorisation of $n$.  
		Using the standard ``{\bf is a divisor of}'' relation, $x | y$ iff $y\ (\md \ x) = 0$, 
		we may write this as $ord_p(n)=\max\{ k : p^k | n\}_{k\in \mathbb N}$.
		
		The order is used to define the {\bf $p$-adic norm} $\| \_ \|_p : \mathbb N \rightarrow \mathbb Q$ by, for all $n\in \mathbb N^+$,
		\[ \| n \|_{p} \ = \ p^{-ord(n)}  \ = \frac{1}{\max\{p^k :  p^k | n \}_{k\in \mathbb N}} \]
		and extended to $\mathbb N$ by taking $\| 0 \|_p = 0$.
		As this is a (non-Archimedean) norm on $\mathbb N$, it defines a  distance, the {\bf $p$-adic distance}. Assuming w.l.o.g. $a\leq b\in \mathbb N$, this is given by  $\| b-a\|_p \ = \  p^{-ord_p(b-a)}$. 
	    \end{definition}
  
    \begin{remark} All the above definitions may be extended to the integers $\mathbb Z$ by taking absolute values, and rationals $\mathbb Q$ since, for all equivalent fractions $\frac{a}{b}=\frac{c}{d}$, we have $\frac{\|a\|_p}{\|b\|_p} =  \frac{\|c\|_p}{\|d\|_p}$.  The $p$-adic reals are then defined as Cauchy sequences w.r.t. the $p$-adic norm on the rationals.
    	\end{remark}
    
    \begin{remark}
      Expanding out the above gives the p-adic distance as 
    \[ \| b-a\|_p\ = \ \frac{1}{ max \{p^k  : a \ (\md \ p^k) = b \} _{k\in \mathbb N}  }  \ \ \ \forall b>a \] 
    Writing it in this form makes the similarity with the {\em residue} and {\em cancellation} of the previous section apparent.
    \end{remark}
    
    We now relate the above elementary definitions to inverse semigroup theory : 
    
    \begin{proposition}\label{polypadic-prop}
    	Recall the embedding $\theta_p: \mathcal P_p \hookrightarrow \mathcal A$. For all $n\in \mathbb N$, 
    	\[  n.\|n\|_p = min \left\{ \theta_p\left(0^k\right)(n)\right\}_{k\in \mathbb N} \]%$, for all $n\in \mathbb N$.\\
    	(Note : we are using $0^k=00\ldots 0$ in the algebraic sense within the free monoid $\{ 0,\ldots,p-1\}^*$, as a string of $k$ copies of $0$). 
    	\end{proposition} % $ $val_p(x) = max \left\{ k\in \mathbb N : x\in dom\left( \theta_p(0)^k \right) \right\}$
    		\begin{proof}
    			By definition, $\theta_p\left(0^k\right)=R_{p,0}^k=R_{p^k,0}$, where 
    			\[ R_{p^k,0}(n) \ = \ \left\{ \begin{array}{lcr} 
    											\frac{n}{p^k} & & \mbox{ when }\  p^k | n \\ & & \\ \bot & & \mbox{otherwise}
    											\end{array}\right.
    											\]
    			Thus $\left\{ \theta_p\left(0^k\right)(n)\right\}_{k\in \mathbb N} \ = \ \left\{ \frac{n}{p^k} : p^k | n \right\}$. 
    			The maximum $k\in \mathbb N$ such that $p^k|n$ gives the minimum of this set; hence $n.\|n\|_p = min \left\{ \theta_p\left(0^k\right)(n)\right\}_{k\in \mathbb N}$, as required.
    			\end{proof}
    		\begin{corollary}
    			The above identity characterises the $p$-adic norm on $\mathbb N$ in terms of the $p$-th polycyclic monoid, as for all $n\neq 0\in \mathbb N$ 
    			\[ \| n \|_p \ = \ min \left\{ \frac{ \theta_p\left(0^k\right)(n)}{n}\right\}_{k\in \mathbb N}  \] 
    		\end{corollary}
    	
    	\begin{remark}
    		One may wonder what is special about the strings 
    		\[ \{ \epsilon, 0 , 00 , 000 , 0000 , \ldots \} \ \subseteq \ \{ 0,\ldots,p-1\}^* \]
    		in the above characterisation of $p$-adic norms?
    		
    		Number-theoretically, the classic theorem of Ostrowski \cite{OS} implies they are very special; algebraically, there is no a priori reason to prefer this subset.  Any infinitary prefix-ordered chain of words from $\{ 0,\ldots,p-1\}^*$ determines a distinct function from $\mathbb N^+$ to $\mathbb R$. This motivates the following definition : 
    	\end{remark}
    
    \begin{definition} Let us denote by $\mathfrak C_p$ the Cantor space of all one-sided infinite words over the set $\{ 0,\ldots , p-1\}$. Given some Cantor point $\Gamma\in \mathfrak C_p$, we define, for all $n\neq 0$ 
    	\[ eval_\Gamma( n) \ = \  min \left\{ \frac{ \theta_p\left(w \right)(n)}{n}\ \ : \ \ w \ \mbox{ is a prefix of }\ \Gamma\in \mathfrak C_p \right\}
    	\]
    	This is always defined, since $\theta_p\left(w \right)(n) = \bot$, for all $len(w)$ sufficiently large, so we are taking the minimum of a finite set in every case.  
    	\end{definition}
    
    	From above, the special case $\Gamma = 00000\ldots $ gives the $p$-adic norm, so $eval_{000\ldots}(n)= \| n\|_p$. The following is straightforward from Proposition \ref{polypadic-prop}, and makes the interpretation of the general case clear :
    	\begin{corollary}
    		Given a non-zero natural $a\in \mathbb N^+$, let us define the Cantor point 
    		$cant(a)\in \mathfrak C_p$ to be the one-sided infinite string
    		\[ cant(a) = a_0a_1\ldots a_x 0000\ldots \]
    		where $n_0n_1\ldots n_x$ is the base-$p$ representation of $n$, and $000\ldots = 0^\omega\in \mathfrak C_p$ is the constantly zero one-sided infinite string. 
    		Then 
    		\[ eval_{cant(a)}(n) = \| n-a\|_p \ \ \forall n>a\in \mathbb N \]
    		\end{corollary} 
    	%	\begin{lemma}
    	%	All {\em negative words} of $P_p$ are strictly contractive, whereas all {\em positive words} are strictly expansive; given arbitrary $u\neq \epsilon\in F_p\hookrightarrow P_p$, then
    	%	\begin{itemize}
    	%		\item $\left\| \theta_p(u)(n)-\theta_p(u)(n') \right\|_p \ < \   \| n-n' \|_p$
    	%		\item $ \left\| \theta_p(u^\ddagger)(n)-\theta_p(u^\ddagger )(n') \right\|_p \ > \   \| n-n' \|_p$
    	%		\end{itemize}
    	%	\end{lemma}
%	\[ \] \[ \] \[ \]
%	\begin{center}{\bf  quite a bit more to write here ...}\end{center}
%	\[ \] \[ \] \[ \]
Thus, considering arbitrary Cantor points leads directly to a curried version of the $p$-adic distance.

Although we could continue and consider questions of convergence, Cauchy sequences, limits, etc., our intention is not to re-construct some well-known theory from inverse semigroup theory! Rather, we wish to highlight an important property, relating to the {\em algebraic} aspects, that  is worthwhile considering.
	
\subsection{An open question}

A great deal of the structure of Nival \& Perot's polycyclic monoids is determined by the fact that they are {\em congruence-free} --- provided $|X|>1$, any monoid homomorphism ${\mathcal P}_X\rightarrow M$ is either an embedding, or maps the whole of ${\mathcal P}_X$ to the identity $1_M\in M$. 

This property  was given a categorical interpretation in \cite{TAC,JHRS,OCL}, where it was related to questions of coherence for associativity, and connected to the well-known fact that Richard Thompson's group $\mathcal F$ has no non-abelian quotients --- which follows for similar categorical reasons.

This raises the following somewhat vaguely stated question :

\begin{center} {\bf What is the number-theoretic significance \\ or interpretation of congruence-freeness,\\  and can we describe this categorically?}  \end{center}

There is a `dual' to the above question, which if anything is even more speculative. A core result for $p$-adic numbers is undoubtedly Ostrowski's theorem on the uniqueness of $p$-adic absolute values. Even more speculatively, we may wonder how this should interpret as inverse semigroup theory?
%(We should, of course, take care to distinguish the algebraic concept of 

\section{Future directions}
Although the algebraic inverse monoid does not appear to have been considered in the literature, there does exist a body of closely related group theory. In \cite{SK}, Stefan Kohl considers a group generated by permutations on $\mathbb N$ that interchange two disjoint\footnote{This disjointness requirement is needed in order to ensure his group is well-defined. This is in contrast to the inverse semigroup theoretic approach, where where the normal form  $R^\ddagger_{c,d}R_{a,b}$ maps monotonically between {\em arbitrary} congruence classes as  $a\mathbb N+b\mapsto c\mathbb N+d$, and is undefined elsewhere.} congruence classes, and act as the identity elsewhere.  This is undoubtedly related to the arithmetic inverse monoid, simply by considering the orthogonal closure of $\mathcal A$ with respect to the natural partial order (see \cite{MVLbool} for the general theory of this technique), then considering the group of units (i.e. globally invertible elements) of the resulting inverse monoid. 

At the very least, this provides a decomposition of the elements of his group into inverse-semigroup theoretic primitives, complete with normal forms and formul\ae\ for composition.  

A more explicitly topological approach is also worth pursuing; the one-sided infinite strings of the Cantor space $\mathfrak C_p$ corresponding to currying the $p$-adic norm on $\mathbb N^+$ are readily identified as the basic clopen sets of the usual topology on $\mathfrak C_p$.  More generally, the idempotents of $\mathcal A$ are the partial identities $\{ Id_{a\mathbb N + b} : b<a\} \cup \{ Id _\emptyset \}$ --- a basis set for the profinite topology on the monoid $(\mathbb N , +)$. The connection between semilattices of idempotents of inverse monoids, and topologies \& locales, is of course well-established.

Finally, the close connection with John Conway's congruential functions \cite{jCon} has not escaped us. Given that a great deal of work on polycyclic monoids has been motivated by their re-discovery as the logicians' `dynamical algebra', in the context of computationally universal systems (\cite{DR,GOI1,GOI2,OCL}) there are undoubtedly some interesting but decidedly non-trivial avenues to explore.

\section*{Acknowledgements} 
Many people have been very helpful, and no work is done in isolation. However, I prefer to add acknowledgements to final versions of papers, in order to take full responsibility for any errors that may occur beforehand.
%I wish to thank the York mathematics algebra/semigroups research group, under the direction of Victoria Gould, for many helpful discussions and seminars. Both Mark Lawson (Heriot-Watt) and Phil Scott (Ottawa) have also made helpful comments and suggestions. 

%I am particularly indebted to Noson Yanofsky (New York) for numerous re-readings of previous drafts, and exceedingly helpful comments on their contents. Many thanks are also due to Matt Brin, for pointing out some relevant work that I had missed in earlier versions.

\bibliographystyle{plain}
\bibliography{arithmetic}

\end{document}